\documentclass{amsart}

\usepackage{graphicx} %incluir figuras

\usepackage{algorithmic, algorithm, xcolor}

\usepackage{blkarray}
\usepackage{cite}
\usepackage{stmaryrd}
\usepackage{hyperref}
\usepackage{amssymb,enumerate}

\renewcommand{\u}{\mathbf{u}}
\renewcommand{\v}{\mathbf{v}}
\renewcommand{\b}{\mathbf{B}}
\newcommand{\gO}{{\color{gray}0}}
\newcommand{\bO}{{\color{blue}1}}
\newcommand{\rO}{{\color{red}1}}
\newcommand{\M}{\mathcal{M}}

\numberwithin{equation}{section}
 \theoremstyle{plain}
\newtheorem{theorem}{Theorem}[section]

\newtheorem{conjecture}[theorem]{Conjecture}

\theoremstyle{remark}

\begin{document}

\makeatletter
\def\imod#1{\allowbreak\mkern10mu({\operator@font mod}\,\,#1)}
\makeatother

\author{Alexander Berkovich}
   \address{Department of Mathematics, University of Florida, 358 Little Hall, Gainesville FL 32611, USA}
   \email{alexb@ufl.edu}
   
\author{Ali Kemal Uncu}
   \address{Austrian Academy of Sciences, Johann Radon Institute for Computational and Applied Mathematics, Altenbergerstrasse 69, 4040 Linz, Austria}
   \email{akuncu@ricam.oeaw.ac.at}

%\scalebox{.9}
\title[Location of the maximum absolute $q$-series coefficient of $(q;q)_n$]{Where do the maximum absolute $q$-series coefficients of $(1-q)(1-q^2)(1-q^3)\dots(1-q^{n-1})(1-q^n)$ occur?}

\begin{abstract} We used the MACH2 supercomputer to study coefficients in the $q$-series expansion of $(1-q)(1-q^2)\dots(1-q^n)$, for all $n\leq 75000$. As a result, we were able to conjecture some periodic properties associated with the before unknown location of the maximum coefficient of these polynomials with odd $n$. Remarkably the observed period is 62{,}624.
\end{abstract}
   
\keywords{Euler Pentagonal Number Theorem, $q$-Pochhammer symbol, Experimental Mathematics, Approximation, Maximum, MACH2, High-Performance Computing}

 \subjclass[2010]{05A15, 05A30, 11Y55, 11Y60, 90C10}
 %%%%	05A15  	Exact enumeration problems, generating functions [See also 33Cxx, 33Dxx]
 % 05A30   	$q$-calculus and related topics
% 11Y60 Computational NT - Evaluation of Constants
%11Y55   	Calculation of integer sequences
%  	90C10   	Integer programming

\thanks{Research of the first author is partly supported by the Simons foundation, Award ID: 308929. Research of the second author is supported by the Austrian Science Fund FWF, SFB50-07, SFB50-09 and SFB50-11 Projects.}

\date{\today}
   \maketitle
   
\section{Introduction}

Let the $q$-Pochhammer symbol be \begin{equation}\label{a_n_i}
(q;q)_n :=\prod_{i= 1}^n (1-q^i) = \sum_{i=0}^{n(n+1)/2} a_{n,i}\, q^i,
\end{equation} for any $n \in \mathbb{Z}_{\geq 0}\cup \{\infty\}$. Euler's famous pentagonal number theorem \cite{Theory_of_Partitions} states that infinite product can be written as a sparse power series with all of its coefficients from the set $\{-1,0,1\}$. More precisely:

\begin{theorem}[Euler's Pentagonal Number Theorem, 1750]\label{EPNT_THM} \begin{align}
\label{EPNT} (q&;q)_\infty = \sum_{n=-\infty} ^{\infty} (-1)^n q^{n(3n-1)/2} = 1 -  q-  q^2+ q^5+  q^7-  q^{12}-  q^{15}+q^{22}+ \dots .
\end{align}
\end{theorem}

In contrast, it is easy to check that the coefficients of the finite products $(q;q)_n$, for $n \in \mathbb{Z}_{\geq 0}$, do not share this property and attain values larger than 1 in size. The authors \cite{BerkovichUncu8} recently proved that the only $(q;q)_n$'s where the series coefficients stay in the set $\{-1,0,1\}$ are when $n=0,1,2,3,$ and $5$. The first maximums of the absolute value of the coefficients are \begin{equation}\label{Heights_list}
1, 1, 1, 1, 2, 1, 2, 2, 2, 2, 3, 2, 4, 3, 3, 4, 6, 5, 6, 7, 8, 8, 10, 11, 16, 16,19, 21, 28, 29, 34,\dots.
\end{equation} 
From now on we will use \textit{maximum absolute coefficient} in place of the wordy ``the maximum of the absolute value of the coefficients" to simplify our language.
 
\begin{figure}[!htb]\caption{Coefficients $a_{250,i}$ plotted on a number line.}\label{Fig1}
\includegraphics[width=\linewidth]{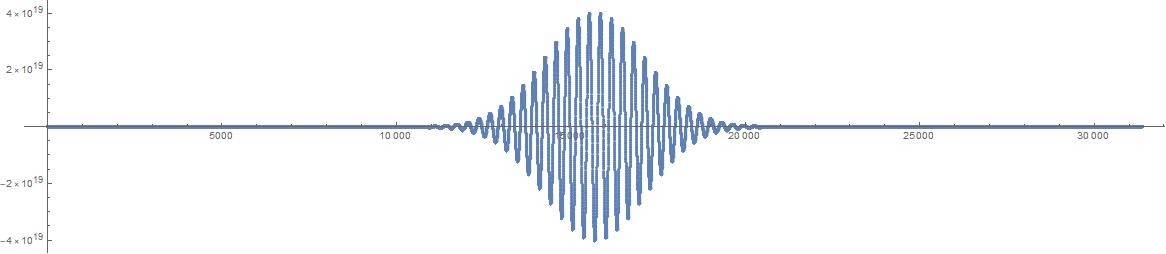}
\end{figure}

In general, when the coefficients $a_{n,i}$ of $(q;q)_n$ are studied, we see tame coefficients at both ends of the polynomial $(q;q)_n$ and a bubble of oscillation involving enormous integers in the middle. We plot the ordered pairs $(i,a_{250,i})$ in Figure~\ref{Fig1}, to show the shape of $(q;q)_{250}$. This is to give an idea of the nature of the coefficients of $(q;q)_n$.
Note that ``tame" coefficients are small only relative to the coefficients that appear at the oscillation. For example, in Figure~\ref{Fig1}, the 5000th coefficient of $(q;q)_{250}$ seems to be close to the $x$-axis, but its value $a_{250,5000}$ is $-7{,}983{,}490$.

Moreover, Sudler \cite{Sudler1, Sudler2} showed that the maximum absolute coefficients' size grows exponentially with $n$: \begin{equation}\label{Sudler_max_size}\M_{n} = e^{K n} + \mathcal{O}( n),\end{equation} where, \[\M_{n}:= \max_{0\leq i\leq n(n+1)/2} |a_{n,i}|.\] Sudler and A. Hurwitz calculated, $K\simeq 0.19861$. On top of that, Wright \cite{Wright} proved that \[e^K = \lim_{n\rightarrow\infty}\frac{\M_n}{ \M_{n-1}}.\]

Although, research on the size of these coefficients existed \cite{Finch, Sudler1, Sudler2}, there has not been an extensive study of the location of the maximum absolute coefficients of $(q;q)_n$ in the literature. When one checks the sequence \eqref{Heights_list} in the Online Encyclopedia of Integer Sequences they are greeted with the open problem~\cite{OEIS} (paraphrased): \begin{quote}\label{quote}
If $n$ is even, then $\M_n$ is the absolute value of the coefficients of $q^{\lfloor n(n+1)/4)\rfloor}$ and $q^{\lceil n(n+1)/4)\rceil}$. If $n$ is odd, it is an open question as for which $i$, $|a_{n,i}|=\M_n$,
\end{quote} where $\lfloor \cdot \rfloor$ and $\lceil \cdot \rceil$ are the classical floor and ceiling functions.

In this work we would like to give an experimental answer to the question regarding the location for the odd $n$. We have carefully calculated and verified the locations of the largest absolute value of the coefficients for all $n$ as well as the maximum absolute coefficients for $n\leq 75 000$. This task had required painstaking attention to details, a lot of patience, and a grand computing power for today's standards. The supercomputer MACH2 was used for this task, which is the largest European installation of its kind with almost 2000 processor cores and 20 terabytes of shared memory. Even with this computing power, the computation was challenging due to the number and the size of the objects that needed to be stored and manipulated at full precision.

\section{Multiplication algorithms used and some related data}

The task of openly multiplying $(q;q)_n$ and locating the maximum of the absolute values of the coefficients started on a reasonably sized RISC server with Maple. One can compute up to a commendable $n=10000$ by multiplying these objects using a modern desktop computer with 8Gb of RAM and a quad core processor of $2.30$Ghz processing speed. In our experience, the Maple gives up around this point. With Mathematica, we were able to take our calculations to around $15000$. This is more or less the upper limit of what computer algebra systems can handle on a normal computer. Once this natural boundary was reached we moved on to programming our own versions of the basic multiplication algorithm; first this was done with Python, later with C++. Python does arbitrary precision arithmetics automatically, and Gnu MultiPrecision Library \cite{GNU} was used for handling large integers in C++.

To save computer memory in calculations, one needs to observe that $a_{n,i}$, as defined in \eqref{a_n_i}, satisfies the relation \begin{equation}\label{symmetry}a_{n,i} = (-1)^n a_{n,n(n+1)/2 - i}.\end{equation} Using \eqref{symmetry}, one can reduce the storage to only half of the coefficients of $(q;q)_n$. Another point is that the multiplication $(1-q^n)(q;q)_{n-1}$ to get $(q;q)_n$ from $(q;q)_{n-1}$ is the same as the addition \begin{equation}\label{multiplication}a_{n,i} = a_{n-1,i} - a_{n-1,i-n},\end{equation}for all $i\geq n$, on the coefficient level. Finally, we note that the first $n-1$ coefficients of $(q;q)_{n-1}$ do not change when multiplied with $(1-q^n)$. Moreover, they coincide with the first $n-1$ coefficients of \eqref{EPNT}. Since the maximum absolute coefficients grow \eqref{Sudler_max_size}, they would not appear in the first $n-1$ coefficients that are $\pm1$ or 0.

The most memory efficient multiplication pseudo-algorithm, that we have come up with, is as follows:

\vspace{.5cm}
\begin{algorithmic}
\STATE Read in the first half of $(q;q)_{Start}$ ($a_{Start,i}$, with $i\leq Start(Start+1)/4$) in an array,
\FOR{$n=Start+1$ \TO $Finish$}
\STATE $Max \leftarrow 0$, $MaxLoc \leftarrow 0$,
\FOR{$i=1$ \TO $\lfloor n/2\rfloor$ } 
\STATE Define and add $a_{n,\lfloor n(n-1)/4\rfloor +i}$ by \eqref{symmetry} and \eqref{multiplication},
\IF{$Max < |a_{n,i}|$} 
\STATE $Max \leftarrow |a_{n,i}|$, $ MaxLoc \leftarrow \lfloor n(n-1)/4\rfloor +i$,
 \ENDIF
\ENDFOR
\FOR{$i = \lfloor n(n-1)/4\rfloor$ \TO $n$}
\STATE $a_{n,i}\leftarrow a_{n-1,i} - a_{n-1,i-n}$,
\IF{$Max < |a_{n,i}|$} 
\STATE $Max \leftarrow |a_{n,i}|$, $ MaxLoc \leftarrow \lfloor n(n-1)/4\rfloor +i$,
 \ENDIF
\ENDFOR
\STATE Record $Max$ and $MaxLoc$,
\ENDFOR
\STATE Record the first half of $(q;q)_{Finish}$ ($a_{Finish,i}$, with $i\leq Finish(Finish+1)/4$).
\end{algorithmic}
\vspace{.3cm}

This algorithm uses a single array of integers and updates this array at each step. Finding the maximum absolute coefficient and its location is done simultaneously with the multiplications. At each step $n$, we add $\lfloor n/2\rfloor$ new integers to the end of the array. In Python, extending arrays can be done automatically, and in C++ a list is used in place of an array. We first calculate the coefficients $a_{n,i}$ for $\lfloor n(n-1)/4\rfloor+1 \leq i\leq \lfloor n(n+1)/4\rfloor $, and calculate all the coefficients up to the half of the degree of the polynomial $(q;q)_n$. This needs to be done using the symmetry of the coefficients \eqref{symmetry}. Then, we start updating the already existing integers $a_{n-1,i}$ with $a_{n,i}$ going backwards from $a_{n,\lfloor n(n-1)/4\rfloor}$, the half of the degree of $(q;q)_{n-1}$. 

One important note here is that an intermediate node $a_{n-1,k}$ for some $n \leq k\leq \lfloor n(n-1)/4\rfloor$ appears in the multiplication of two terms: $a_{n,k}$ and $a_{n,k+n}$. Hence, this node should not be updated before these two multiplications are handled. This throws a ratchet in the gears of parallelization of this computation. Although this algorithm is memory efficient (as much as it can be), it is not programmed with parallel computations in mind. We can parallelize this algorithm by splitting the multiplications into residue classes modulo $n$ for each multiplication, but since the number of residue classes changes at each step, this possibility did not look beneficial.

The above mentioned memory efficient algorithm is how we carried the calculations towards $(q;q)_{50000}$. With $1.5$ Terabytes of RAM, the largest server at RISC, qftquad8, was able to save an image of $(q;q)_{45000}$ for us to carry our calculations to MACH2. It also recorded the maximum absolute coefficients and their locations till about 49000. These calculations also allowed us to check and confirm their validity when we switched to the new algorithm for multiplication we used in MACH2. The supercomputer MACH2 is a shared computer with 1728 cores and 20 terabytes of accessible memory. Therefore, for speeding up the calculations we moved onto a version of the above described multiplication algorithm that lavishly uses the available memory.

\vspace{.3cm}
\begin{algorithmic}
\STATE Initialize two arrays $Arr1$ and $Arr2$,
\STATE Read in the first half of $(q;q)_{Start}$ in $Arr1$, and the first $Start$ amount of entries in $Arr2$,
\FOR{$n=Start+1$ \TO $Finish$}
\STATE $Max \leftarrow 0$, $MaxLoc \leftarrow 0$,
\FOR{$i=1$ \TO $n$ } 
\STATE In parallel, define and add $a_{n,\lfloor n(n-1)/4\rfloor +i}$ by \eqref{symmetry} and \eqref{multiplication} to $Arr2$,
\IF{$Max < |a_{n,i}|$} 
\STATE $Max \leftarrow |a_{n,i}|$, $ MaxLoc \leftarrow \lfloor n(n-1)/4\rfloor +i$,
 \ENDIF
\ENDFOR
\FOR{$i = n$ \TO $\lfloor n(n-1)/4\rfloor$}
\STATE In Parallel, $a_{n,i}\leftarrow a_{n-1,i} - a_{n-1,i-n}$ in $Arr2$,
\IF{$Max < |a_{n,i}|$} 
\STATE $Max \leftarrow |a_{n,i}|$, $ MaxLoc \leftarrow \lfloor n(n-1)/4\rfloor +i$,
 \ENDIF
\ENDFOR
\STATE Record $Max$ and $MaxLoc$,
\STATE Reverse the roles of $Arr1$ and $Arr2$, and repeat the multiplications,
\ENDFOR
\STATE Record the first half of $(q;q)_{Finish}$.
\end{algorithmic}
\vspace{.3cm}

Since the reading is done in one array and the writing is done in another, the multiplication can easily be parallelized. The price paid is the double the amount of memory needed to keep two separate lists of huge integers. Moreover, to speed up the process ever so slightly, we can use arrays of the final needed sizes for $Arr1$ and $Arr2$ instead of lists. This small change reduces the complexities that pointers of C++ introduce greatly.

The calculations in the supercomputer MACH2 have been carried over using a 10 terabyte capsule with hundreds of cores working in parallel. In this phase of the calculation, the first half of the $(q;q)_n$'s were recorded in increments of 5000 from 45000 till 75000. Some specifics of the MACH2 Supercomputer are as follows: CPU Type: Dodeka-Core Haswell CPU "Intel(R) Xeon(R) CPU E5-4650 v3 @ 2.10GHz"; Linux Kernel: 4.4.120; Platform: x86$\underline{\text{ }}$64; GNU libc 2.22; GNU Multiprecision Library Version: GMP 5.1.3; GNU Compiler Collection: 4.8.5\footnote{See \url{https://www.risc.jku.at/projects/mach2}.}. 

Towards the of the calculations, the $100$ consecutive multiplications with identification of the maximum absolute coefficient took around 1700 Mega-GMP-Operations per second computing time, on average. The actual calculations from start to finish with the recording of the checkpoint polynomials took around 6 months. With interruptions, designing/coding challenges, data migrations, and checks, this project took the authors around two full years.

\section{The location of the maximum absolute coefficients}

Let $L(n)$ denote the lowest exponent $i$ of $q$ in the $(q;q)_n$ polynomial, where the coefficient of $q^i$ is the maximum absolute coefficient $\M_n$. When explicitly expanding the polynomials $(q;q)_n$, the maximum of the absolute value of the coefficients starts to appear in one or two places consistently for $n\geq 34$. We would like to start by acknowledging that, for all $n\geq 17$, $L(2n) = \lfloor n(n+1)/4)\rfloor$ or, in other words, the maximum absolute coefficients of $(q;q)_{2n}$ appear as the coefficient(s) of $q^{\lfloor n(n+1)/4)\rfloor}$ and $q^{\lceil n(n+1)/4)\rceil}$. Now with our explicit computations, this claim (which appears in the Online Encyclopedia of Integer Sequences A160089) extends to all even~$n \leq 75000$.  

For 0 modulo 4 cases the maximum absolute coefficients are the absolute maximum, and for 2 modulo 4 cases the maximum absolute coefficients are the negative of the absolute minimum of $(q;q)_n$. For the odd $n$, both the maximum of the absolute value of the coefficients and its negative appear in $(q;q)_n$ due to \eqref{symmetry}. Furthermore, we make the claim: 
\begin{conjecture} For $n\geq 35$, if $n\equiv 1\text{ (mod 4)}$, absolute maximum occurs only once as the coefficient of some $q^i$, where $i< n(n+1)/4$ and if~$n\equiv 3\text{ (mod 4)}$ the absolute maximum occurs only once as the coefficient of some $q^i$, where $i> n(n+1)/4$.
\end{conjecture}
There are some initial cases where the maximum absolute coefficient appears as the polynomial coefficients of $(q;q)_n$ more than twice. For example, in the $n=33$ case the absolute maximum coefficient $56$ appears as the coefficients of $q^{270}$ and $q^{272}$ and the absolute minimum coefficient $-56$ appears as the coefficients of $q^{289}$ and $q^{291}$, so the maximum absolute coefficient appears a total of 4 times as a coefficient of the polynomial $(q;q)_{33}$. We see that for $n\geq 35$, the maximum absolute coefficient appears as coefficients of only two terms in the $q$-series expansion of the rising factorial.

We will split the cases for the odd values of $n$ in two groups modulo 4. Before doing so, we define the $D_n$ as the half difference between the two locations of the maximum absolute coefficients for a given odd $n\geq 35$. This is the same as saying how far off the maximum absolute coefficient is from the half degree $n(n+1)/4$ of the polynomial $(q;q)_n$. Hence, knowing this value and $n$ is enough to recover the location of the maximum absolute coefficients. Also we define the canonical difference \[E_n:= D_{n}-D_{n-4},\] for all odd $n\geq 35$. 

Let $n\geq 35$ be an odd integer, and assume that all necessary $E_{i}$ and a $D_m$ value for some $m\leq n$, with $n\equiv m$ modulo 4 are given. Then, one can recover the location of the maximum absolute coefficients for $(q;q)_n$ as \begin{equation}\label{Location} L(n)=\frac{n(n+1)}{4} - \left(D_m + \sum_{i=1}^{(n-m)/4} E_{m+4i}\right),
\end{equation} where $a_{n,i}$ is defined as in \eqref{a_n_i}.

The first observation we have made and confirmed for $n \leq 75000$ is:

\begin{conjecture}\label{e_size} For all odd $n\geq 35$, $E_n\in\{0,1,2\}$. Furthermore, $E_n >0$ for all odd $n\geq 61$.
\end{conjecture}

By looking at the values $n\geq 57$ and $\geq 87$ for $1$ and $3$ modulo 4, we list $E_n$'s. The patterns starting from these beginning points look highly periodic with period 19, although this is sadly not the case. The 19 length patterns of 1's and 2's that $E_n$'s change slightly as the calculations are carried. We match the 19 length patterns that the calculations yield with letters and form an alphabet with 20 letters. The letters that we see in the calculations and their 19 consecutive $E_n$ value equivalents are given in Table~\ref{Table_1}.

\begin{table}[!htb]
\caption{The selected alphabet associated with the $E_m$ values.}\label{Table_1}\vspace{-.5cm}
{\footnotesize\begin{align*}
a &= \overbrace{2112111211121112112}^{19}, \hspace{-2.5cm} &k= \overbrace{1211121112111211121}^{19},\\
b &= 1112111211121112112, \hspace{-2.5cm} &l= 1211121112111211211,\\
c &= 1112111211121121112, \hspace{-2.5cm} &m= 1211121112112111211,\\
d &= 1112111211211121112, \hspace{-2.5cm} &n= 1211121121112111211,\\
e &= 1112112111211121112, \hspace{-2.5cm} &o= 1211211121112111211,\\
f &= 1121112111211121112, \hspace{-2.5cm} &p= 2111211121112111211,\\
g &= 1121112111211121121, \hspace{-2.5cm} &q= 2111211121112112111,\\
h &= 1121112111211211121, \hspace{-2.5cm} &r= 2111211121121112111,\\
i &= 1121112112111211121, \hspace{-2.5cm} &s= 2111211211121112111,\\
j &= 1121121112111211121, \hspace{-2.5cm} &t= 2112111211121112111.\\
\end{align*}}\vspace{-1.1cm}
\end{table}

Note that each letter represents 19 consecutive $E_n$ values (1 or 2) written together without commas to keep the notation simple. The calculation of any letter requires the expansion of $19\times 4 = 76$ consecutive $(q;q)_n$ polynomials and finding the location of the maximum absolute coefficient and its location of each odd index $n$. This alphabet makes it possible to find the individual $E_n$ values. Furthermore, together with~\eqref{Location}, this makes it possible to find the locations of the maximum absolute coefficients. Also note that $a$ and $k$ are symmetric with respect to their centers and $a$ is also special that it is the only letter with six 2 values, whereas all the other words come with five 2 values. It is easy to see that the letters $b,c,d,e,f,g,h,i,$ and $j$ read from right to left are the letters $t,s,r,q,p,o,n,m,$ and $l$, respectively.

These letters are observed to be generating almost periodic words starting from $n\geq 209$ and $\geq 391$  for $1$ and $3$ modulo 4 values, respectively. We will write exponents of the letters defined in Table~\ref{Table_1} to indicate the number of times a letter has occurred consecutively. From the said starting point we see that the 1 modulo 4 differences give rise to the string \[\underbrace{a^1 b^4 c^4 d^4 e^4 f^3 g^4 h^5 i^4 j^4 k^3 l^4 m^4 n^4 o^5 p^3 q^4 r^4 s^4 t^3}_{\text{word 1}} \underbrace{a^1 b^4 c^4 d^4 e^4 f^3 g^4 h^4 i^5 j^4 k^3 \dots}_{\text{word 2, etc.}}.\] 

We see that the letter $a$ appears a single time in each word. This phenomenon is the same for the 3 modulo 4 case. We break this string into words for $1$ and $3$ mod 4 to start every line with an $a$ and write the words in Table~\ref{Table_words}; 1 and 3 mod 4 cases on the left side and the right side, respectively. In Table~\ref{Table_words}, we also include a line to break period blocks of the size $11\times 20$ box both odd residue classes mod 4.

\begin{table}[h]
\caption{The words for $n\equiv1$ and 3 mod 4 cases associated with the alphabet of Table~\ref{Table_1}.}\label{Table_words}{\tiny
\begin{tabular}{c|c}
$\begin{array}{l}
\begin{array}{llllllllllllllllllll}
a^1 &\hspace{-3mm} b^4 &\hspace{-3mm} c^4 &\hspace{-3mm} d^4 &\hspace{-3mm} e^4 &\hspace{-3mm} f^3 &\hspace{-3mm} g^4 &\hspace{-3mm} h^5 &\hspace{-3mm} i^4 &\hspace{-3mm} j^4 &\hspace{-3mm} k^3 &\hspace{-3mm} l^4 &\hspace{-3mm} m^4 &\hspace{-3mm} n^4 &\hspace{-3mm} o^5 &\hspace{-3mm} p^3 &\hspace{-3mm} q^4 &\hspace{-3mm} r^4 &\hspace{-3mm} s^4 &\hspace{-3mm} t^3 \\ \hline
\end{array}\\
\hspace{-5.2mm}11\left\{\begin{array}{llllllllllllllllllll}
a^1 &\hspace{-3mm} b^4 &\hspace{-3mm} c^4 &\hspace{-3mm} d^4 &\hspace{-3mm} e^4 &\hspace{-3mm} f^3 &\hspace{-3mm} g^4 &\hspace{-3mm} h^4 &\hspace{-3mm} i^5 &\hspace{-3mm} j^4 &\hspace{-3mm} k^3 &\hspace{-3mm} l^4 &\hspace{-3mm} m^4 &\hspace{-3mm} n^4 &\hspace{-3mm} o^5 &\hspace{-3mm} p^3 &\hspace{-3mm} q^4 &\hspace{-3mm} r^4 &\hspace{-3mm} s^4 &\hspace{-3mm} t^3 \\
a^1 &\hspace{-3mm} b^3 &\hspace{-3mm} c^5 &\hspace{-3mm} d^4 &\hspace{-3mm} e^4 &\hspace{-3mm} f^3 &\hspace{-3mm} g^4 &\hspace{-3mm} h^4 &\hspace{-3mm} i^5 &\hspace{-3mm} j^4 &\hspace{-3mm} k^3 &\hspace{-3mm} l^4 &\hspace{-3mm} m^4 &\hspace{-3mm} n^4 &\hspace{-3mm} o^4 &\hspace{-3mm} p^4 &\hspace{-3mm} q^4 &\hspace{-3mm} r^4 &\hspace{-3mm} s^4 &\hspace{-3mm} t^3 \\
a^1 &\hspace{-3mm} b^3 &\hspace{-3mm} c^5 &\hspace{-3mm} d^4 &\hspace{-3mm} e^4 &\hspace{-3mm} f^3 &\hspace{-3mm} g^4 &\hspace{-3mm} h^4 &\hspace{-3mm} i^4 &\hspace{-3mm} j^5 &\hspace{-3mm} k^3 &\hspace{-3mm} l^4 &\hspace{-3mm} m^4 &\hspace{-3mm} n^4 &\hspace{-3mm} o^4 &\hspace{-3mm} p^4 &\hspace{-3mm} q^4 &\hspace{-3mm} r^4 &\hspace{-3mm} s^4 &\hspace{-3mm} t^3 \\
a^1 &\hspace{-3mm} b^3 &\hspace{-3mm} c^4 &\hspace{-3mm} d^5 &\hspace{-3mm} e^4 &\hspace{-3mm} f^3 &\hspace{-3mm} g^4 &\hspace{-3mm} h^4 &\hspace{-3mm} i^4 &\hspace{-3mm} j^5 &\hspace{-3mm} k^3 &\hspace{-3mm} l^4 &\hspace{-3mm} m^4 &\hspace{-3mm} n^4 &\hspace{-3mm} o^4 &\hspace{-3mm} p^3 &\hspace{-3mm} q^5 &\hspace{-3mm} r^4 &\hspace{-3mm} s^4 &\hspace{-3mm} t^3 \\
a^1 &\hspace{-3mm} b^3 &\hspace{-3mm} c^4 &\hspace{-3mm} d^5 &\hspace{-3mm} e^4 &\hspace{-3mm} f^3 &\hspace{-3mm} g^4 &\hspace{-3mm} h^4 &\hspace{-3mm} i^4 &\hspace{-3mm} j^4 &\hspace{-3mm} k^4 &\hspace{-3mm} l^4 &\hspace{-3mm} m^4 &\hspace{-3mm} n^4 &\hspace{-3mm} o^4 &\hspace{-3mm} p^3 &\hspace{-3mm} q^4 &\hspace{-3mm} r^5 &\hspace{-3mm} s^4 &\hspace{-3mm} t^3 \\
a^1 &\hspace{-3mm} b^3 &\hspace{-3mm} c^4 &\hspace{-3mm} d^4 &\hspace{-3mm} e^5 &\hspace{-3mm} f^3 &\hspace{-3mm} g^4 &\hspace{-3mm} h^4 &\hspace{-3mm} i^4 &\hspace{-3mm} j^4 &\hspace{-3mm} k^3 &\hspace{-3mm} l^5 &\hspace{-3mm} m^4 &\hspace{-3mm} n^4 &\hspace{-3mm} o^4 &\hspace{-3mm} p^3 &\hspace{-3mm} q^4 &\hspace{-3mm} r^5 &\hspace{-3mm} s^4 &\hspace{-3mm} t^3 \\
a^1 &\hspace{-3mm} b^3 &\hspace{-3mm} c^4 &\hspace{-3mm} d^4 &\hspace{-3mm} e^4 &\hspace{-3mm} f^4 &\hspace{-3mm} g^4 &\hspace{-3mm} h^4 &\hspace{-3mm} i^4 &\hspace{-3mm} j^4 &\hspace{-3mm} k^3 &\hspace{-3mm} l^5 &\hspace{-3mm} m^4 &\hspace{-3mm} n^4 &\hspace{-3mm} o^4 &\hspace{-3mm} p^3 &\hspace{-3mm} q^4 &\hspace{-3mm} r^4 &\hspace{-3mm} s^5 &\hspace{-3mm} t^3 \\
a^1 &\hspace{-3mm} b^3 &\hspace{-3mm} c^4 &\hspace{-3mm} d^4 &\hspace{-3mm} e^4 &\hspace{-3mm} f^4 &\hspace{-3mm} g^4 &\hspace{-3mm} h^4 &\hspace{-3mm} i^4 &\hspace{-3mm} j^4 &\hspace{-3mm} k^3 &\hspace{-3mm} l^4 &\hspace{-3mm} m^5 &\hspace{-3mm} n^4 &\hspace{-3mm} o^4 &\hspace{-3mm} p^3 &\hspace{-3mm} q^4 &\hspace{-3mm} r^4 &\hspace{-3mm} s^5 &\hspace{-3mm} t^3 \\
a^1 &\hspace{-3mm} b^3 &\hspace{-3mm} c^4 &\hspace{-3mm} d^4 &\hspace{-3mm} e^4 &\hspace{-3mm} f^3 &\hspace{-3mm} g^5 &\hspace{-3mm} h^4 &\hspace{-3mm} i^4 &\hspace{-3mm} j^4 &\hspace{-3mm} k^3 &\hspace{-3mm} l^4 &\hspace{-3mm} m^5 &\hspace{-3mm} n^4 &\hspace{-3mm} o^4 &\hspace{-3mm} p^3 &\hspace{-3mm} q^4 &\hspace{-3mm} r^4 &\hspace{-3mm} s^4 &\hspace{-3mm} t^4 \\
a^1 &\hspace{-3mm} b^3 &\hspace{-3mm} c^4 &\hspace{-3mm} d^4 &\hspace{-3mm} e^4 &\hspace{-3mm} f^3 &\hspace{-3mm} g^5 &\hspace{-3mm} h^4 &\hspace{-3mm} i^4 &\hspace{-3mm} j^4 &\hspace{-3mm} k^3 &\hspace{-3mm} l^4 &\hspace{-3mm} m^4 &\hspace{-3mm} n^5 &\hspace{-3mm} o^4 &\hspace{-3mm} p^3 &\hspace{-3mm} q^4 &\hspace{-3mm} r^4 &\hspace{-3mm} s^4 &\hspace{-3mm} t^4 \\
a^1 &\hspace{-3mm} b^3 &\hspace{-3mm} c^4 &\hspace{-3mm} d^4 &\hspace{-3mm} e^4 &\hspace{-3mm} f^3 &\hspace{-3mm} g^4 &\hspace{-3mm} h^5 &\hspace{-3mm} i^4 &\hspace{-3mm} j^4 &\hspace{-3mm} k^3 &\hspace{-3mm} l^4 &\hspace{-3mm} m^4 &\hspace{-3mm} n^5 &\hspace{-3mm} o^4 &\hspace{-3mm} p^3 &\hspace{-3mm} q^4 &\hspace{-3mm} r^4 &\hspace{-3mm} s^4 &\hspace{-3mm} t^3 \\
\end{array}\right.\\
\begin{array}{llllllllllllllllllll} \hline
a^1 &\hspace{-3mm} b^4 &\hspace{-3mm} c^4 &\hspace{-3mm} d^4 &\hspace{-3mm} e^4 &\hspace{-3mm} f^3 &\hspace{-3mm} g^4 &\hspace{-3mm} h^4 &\hspace{-3mm} i^5 &\hspace{-3mm} j^4 &\hspace{-3mm} k^3 &\hspace{-3mm} l^4 &\hspace{-3mm} m^4 &\hspace{-3mm} n^4 &\hspace{-3mm} o^5 &\hspace{-3mm} p^3 &\hspace{-3mm} q^4 &\hspace{-3mm} r^4 &\hspace{-3mm} s^4 &\hspace{-3mm} t^3 \\
a^1 &\hspace{-3mm} b^3 &\hspace{-3mm} c^5 &\hspace{-4mm} \dots 
\end{array}
\end{array}$
&\hspace{-.1cm}
$\begin{array}{l}
\hspace{-.5mm}\left.\begin{array}{llllllllllllllllllll}
a^1 &\hspace{-3mm} b^3 &\hspace{-3mm} c^4 &\hspace{-3mm} d^5 &\hspace{-3mm} e^4 &\hspace{-3mm} f^3 &\hspace{-3mm} g^4 &\hspace{-3mm} h^4 &\hspace{-3mm} i^4 &\hspace{-3mm} j^4 &\hspace{-3mm} k^4 &\hspace{-3mm} l^4 &\hspace{-3mm} m^4 &\hspace{-3mm} n^4 &\hspace{-3mm} o^4 &\hspace{-3mm} p^3 &\hspace{-3mm} q^5 &\hspace{-3mm} r^4 &\hspace{-3mm} s^4 &\hspace{-3mm} t^3 \\
a^1 &\hspace{-3mm} b^3 &\hspace{-3mm} c^4 &\hspace{-3mm} d^4 &\hspace{-3mm} e^5 &\hspace{-3mm} f^3 &\hspace{-3mm} g^4 &\hspace{-3mm} h^4 &\hspace{-3mm} i^4 &\hspace{-3mm} j^4 &\hspace{-3mm} k^4 &\hspace{-3mm} l^4 &\hspace{-3mm} m^4 &\hspace{-3mm} n^4 &\hspace{-3mm} o^4 &\hspace{-3mm} p^3 &\hspace{-3mm} q^4 &\hspace{-3mm} r^5 &\hspace{-3mm} s^4 &\hspace{-3mm} t^3 \\
a^1 &\hspace{-3mm} b^3 &\hspace{-3mm} c^4 &\hspace{-3mm} d^4 &\hspace{-3mm} e^5 &\hspace{-3mm} f^3 &\hspace{-3mm} g^4 &\hspace{-3mm} h^4 &\hspace{-3mm} i^4 &\hspace{-3mm} j^4 &\hspace{-3mm} k^3 &\hspace{-3mm} l^5 &\hspace{-3mm} m^4 &\hspace{-3mm} n^4 &\hspace{-3mm} o^4 &\hspace{-3mm} p^3 &\hspace{-3mm} q^4 &\hspace{-3mm} r^4 &\hspace{-3mm} s^5 &\hspace{-3mm} t^3 \\
a^1 &\hspace{-3mm} b^3 &\hspace{-3mm} c^4 &\hspace{-3mm} d^4 &\hspace{-3mm} e^4 &\hspace{-3mm} f^4 &\hspace{-3mm} g^4 &\hspace{-3mm} h^4 &\hspace{-3mm} i^4 &\hspace{-3mm} j^4 &\hspace{-3mm} k^3 &\hspace{-3mm} l^4 &\hspace{-3mm} m^5 &\hspace{-3mm} n^4 &\hspace{-3mm} o^4 &\hspace{-3mm} p^3 &\hspace{-3mm} q^4 &\hspace{-3mm} r^4 &\hspace{-3mm} s^5 &\hspace{-3mm} t^3 \\
a^1 &\hspace{-3mm} b^3 &\hspace{-3mm} c^4 &\hspace{-3mm} d^4 &\hspace{-3mm} e^4 &\hspace{-3mm} f^3 &\hspace{-3mm} g^5 &\hspace{-3mm} h^4 &\hspace{-3mm} i^4 &\hspace{-3mm} j^4 &\hspace{-3mm} k^3 &\hspace{-3mm} l^4 &\hspace{-3mm} m^5 &\hspace{-3mm} n^4 &\hspace{-3mm} o^4 &\hspace{-3mm} p^3 &\hspace{-3mm} q^4 &\hspace{-3mm} r^4 &\hspace{-3mm} s^4 &\hspace{-3mm} t^4 \\
a^1 &\hspace{-3mm} b^3 &\hspace{-3mm} c^4 &\hspace{-3mm} d^4 &\hspace{-3mm} e^4 &\hspace{-3mm} f^3 &\hspace{-3mm} g^5 &\hspace{-3mm} h^4 &\hspace{-3mm} i^4 &\hspace{-3mm} j^4 &\hspace{-3mm} k^3 &\hspace{-3mm} l^4 &\hspace{-3mm} m^4 &\hspace{-3mm} n^5 &\hspace{-3mm} o^4 &\hspace{-3mm} p^3 &\hspace{-3mm} q^4 &\hspace{-3mm} r^4 &\hspace{-3mm} s^4 &\hspace{-3mm} t^4 \\
a^1 &\hspace{-3mm} b^3 &\hspace{-3mm} c^4 &\hspace{-3mm} d^4 &\hspace{-3mm} e^4 &\hspace{-3mm} f^3 &\hspace{-3mm} g^4 &\hspace{-3mm} h^5 &\hspace{-3mm} i^4 &\hspace{-3mm} j^4 &\hspace{-3mm} k^3 &\hspace{-3mm} l^4 &\hspace{-3mm} m^4 &\hspace{-3mm} n^5 &\hspace{-3mm} o^4 &\hspace{-3mm} p^3 &\hspace{-3mm} q^4 &\hspace{-3mm} r^4 &\hspace{-3mm} s^4 &\hspace{-3mm} t^3 \\
a^1 &\hspace{-3mm} b^4 &\hspace{-3mm} c^4 &\hspace{-3mm} d^4 &\hspace{-3mm} e^4 &\hspace{-3mm} f^3 &\hspace{-3mm} g^4 &\hspace{-3mm} h^5 &\hspace{-3mm} i^4 &\hspace{-3mm} j^4 &\hspace{-3mm} k^3 &\hspace{-3mm} l^4 &\hspace{-3mm} m^4 &\hspace{-3mm} n^4 &\hspace{-3mm} o^5 &\hspace{-3mm} p^3 &\hspace{-3mm} q^4 &\hspace{-3mm} r^4 &\hspace{-3mm} s^4 &\hspace{-3mm} t^3 \\
a^1 &\hspace{-3mm} b^4 &\hspace{-3mm} c^4 &\hspace{-3mm} d^4 &\hspace{-3mm} e^4 &\hspace{-3mm} f^3 &\hspace{-3mm} g^4 &\hspace{-3mm} h^4 &\hspace{-3mm} i^5 &\hspace{-3mm} j^4 &\hspace{-3mm} k^3 &\hspace{-3mm} l^4 &\hspace{-3mm} m^4 &\hspace{-3mm} n^4 &\hspace{-3mm} o^4 &\hspace{-3mm} p^4 &\hspace{-3mm} q^4 &\hspace{-3mm} r^4 &\hspace{-3mm} s^4 &\hspace{-3mm} t^3 \\
a^1 &\hspace{-3mm} b^3 &\hspace{-3mm} c^5 &\hspace{-3mm} d^4 &\hspace{-3mm} e^4 &\hspace{-3mm} f^3 &\hspace{-3mm} g^4 &\hspace{-3mm} h^4 &\hspace{-3mm} i^4 &\hspace{-3mm} j^5 &\hspace{-3mm} k^3 &\hspace{-3mm} l^4 &\hspace{-3mm} m^4 &\hspace{-3mm} n^4 &\hspace{-3mm} o^4 &\hspace{-3mm} p^4 &\hspace{-3mm} q^4 &\hspace{-3mm} r^4 &\hspace{-3mm} s^4 &\hspace{-3mm} t^3 \\
a^1 &\hspace{-3mm} b^3 &\hspace{-3mm} c^4 &\hspace{-3mm} d^5 &\hspace{-3mm} e^4 &\hspace{-3mm} f^3 &\hspace{-3mm} g^4 &\hspace{-3mm} h^4 &\hspace{-3mm} i^4 &\hspace{-3mm} j^5 &\hspace{-3mm} k^3 &\hspace{-3mm} l^4 &\hspace{-3mm} m^4 &\hspace{-3mm} n^4 &\hspace{-3mm} o^4 &\hspace{-3mm} p^3 &\hspace{-3mm} q^5 &\hspace{-3mm} r^4 &\hspace{-3mm} s^4 &\hspace{-3mm} t^3
\end{array} \right\}11\\
\begin{array}{llllllllllllllllllll}\hline
a^1 &\hspace{-3mm} b^3 &\hspace{-3mm} c^4 &\hspace{-3mm} d^5 &\hspace{-3mm} e^4 &\hspace{-3mm} f^3 &\hspace{-3mm} g^4 &\hspace{-3mm} h^4 &\hspace{-3mm} i^4 &\hspace{-3mm} j^4 &\hspace{-3mm} k^4 &\hspace{-3mm} l^4 &\hspace{-3mm} m^4 &\hspace{-3mm} n^4 &\hspace{-3mm} o^4 &\hspace{-3mm} p^3 &\hspace{-3mm} q^5 &\hspace{-3mm} r^4 &\hspace{-3mm} s^4 &\hspace{-3mm} t^3 \\
a^1 &\hspace{-3mm} b^3 &\hspace{-3mm} c^4 &\hspace{-3mm} d^4 &\hspace{-3mm} e^5 &\hspace{-3mm} f^3 &\hspace{-3mm} g^4 &\hspace{-3mm} h^4 &\hspace{-3mm} i^4 &\hspace{-3mm} j^4 &\hspace{-3mm} k^4 &\hspace{-3mm} l^4 &\hspace{-3mm} m^4 &\hspace{-3mm} n^4 &\hspace{-3mm} o^4 &\hspace{-3mm} p^3 &\hspace{-3mm} q^4 &\hspace{-3mm} r^5 &\hspace{-3mm} s^4 &\hspace{-3mm} t^3 \\
a^1 &\hspace{-3mm} b^3 &\hspace{-3mm} c^4 &\hspace{-4mm} \dots 
\end{array}
\end{array}$
\end{tabular}}
\end{table}

What one can observe from the encoding of Table~\ref{Table_words} is that there is the base word \begin{equation}\label{baseword}a^1\ b^3\ c^4\ d^4\ e^4\ f^3\ g^4\ h^4\ i^4\ j^4\ k^3\ l^4\ m^4\ n^4\ o^4\ p^3\ q^4\ r^4\ s^4\ t^3 \end{equation} in each line and an overlying, moving perturbation. It is observed that the first line of the 1 modulo 4 case does not fit any period block and stays as an outlier. We interpret this as the asymptotic behaviour of the perturbation kicking in action a little later than the 3 modulo 4 case. More importantly, this perturbation is observably periodic. To express the observed structure better, we represent the base word \eqref{baseword} as a vector, where we write the frequencies of the letters as a vector of 20 entries \begin{equation}\label{basevector}
\b=(1,3,4,4,4,3,4,4,4,4,3,4,4,4,4,3,4,4,4,3).
\end{equation}
Then Table~\ref{Table_words} can be expressed as adding some perturbation vectors to $\b$. We define the perturbation vectors for the 1 and 3 modulo 4 classes side-by-side in Table~\ref{TablE_perturbation_vectors}.

\begin{table}[!htb]
\caption{The perturbation vectors of a full observed period for 1 and 3 modulo 4 cases side-by-side.}\label{TablE_perturbation_vectors}\vspace{-.5cm}
{\footnotesize\[\begin{array}{l|l}
\u_1  \:\,= (\gO,\rO,\gO,\gO,\gO,\gO,\gO,\gO,\bO,\gO,\gO,\gO,\gO,\gO,\rO,\gO,\gO,\gO,\gO,\gO),  	&\v_1\:\,= (\gO,\gO,\gO,\bO,\gO,\gO,\gO,\gO,\gO,\gO,\bO,\gO,\gO,\gO,\gO,\gO,\bO,\gO,\gO,\gO),\\
\u_2 \:\,= (\gO,\gO,\bO,\gO,\gO,\gO,\gO,\gO,\bO,\gO,\gO,\gO,\gO,\gO,\gO,\bO,\gO,\gO,\gO,\gO),  	&\v_2\:\,= (\gO,\gO,\gO,\gO,\bO,\gO,\gO,\gO,\gO,\gO,\bO,\gO,\gO,\gO,\gO,\gO,\gO,\rO,\gO,\gO),\\
\u_3 \:\,= (\gO,\gO,\bO,\gO,\gO,\gO,\gO,\gO,\gO,\bO,\gO,\gO,\gO,\gO,\gO,\bO,\gO,\gO,\gO,\gO),  	&\v_3\:\,= (\gO,\gO,\gO,\gO,\bO,\gO,\gO,\gO,\gO,\gO,\gO,\rO,\gO,\gO,\gO,\gO,\gO,\gO,\bO,\gO),\\
\u_4 \:\,= (\gO,\gO,\gO,\bO,\gO,\gO,\gO,\gO,\gO,\bO,\gO,\gO,\gO,\gO,\gO,\gO,\rO,\gO,\gO,\gO), 	&\v_4\:\,= (\gO,\gO,\gO,\gO,\gO,\rO,\gO,\gO,\gO,\gO,\gO,\gO,\bO,\gO,\gO,\gO,\gO,\gO,\bO,\gO),\\
\u_5 \:\,= (\gO,\gO,\gO,\bO,\gO,\gO,\gO,\gO,\gO,\gO,\rO,\gO,\gO,\gO,\gO,\gO,\gO,\bO,\gO,\gO),  	&\v_5\:\,= (\gO,\gO,\gO,\gO,\gO,\gO,\bO,\gO,\gO,\gO,\gO,\gO,\bO,\gO,\gO,\gO,\gO,\gO,\gO,\bO),\\
\u_6 \:\,= (\gO,\gO,\gO,\gO,\rO,\gO,\gO,\gO,\gO,\gO,\gO,\bO,\gO,\gO,\gO,\gO,\gO,\bO,\gO,\gO),  	&\v_6\:\,= (\gO,\gO,\gO,\gO,\gO,\gO,\bO,\gO,\gO,\gO,\gO,\gO,\gO,\bO,\gO,\gO,\gO,\gO,\gO,\bO),\\
\u_7 \:\,= (\gO,\gO,\gO,\gO,\gO,\bO,\gO,\gO,\gO,\gO,\gO,\bO,\gO,\gO,\gO,\gO,\gO,\gO,\bO,\gO),  	&\v_7\:\,= (\gO,\gO,\gO,\gO,\gO,\gO,\gO,\bO,\gO,\gO,\gO,\gO,\gO,\bO,\gO,\gO,\gO,\gO,\gO,\gO),\\
\u_8 \:\,= (\gO,\gO,\gO,\gO,\gO,\bO,\gO,\gO,\gO,\gO,\gO,\gO,\bO,\gO,\gO,\gO,\gO,\gO,\bO,\gO),  	&\v_8\:\,= (\gO,\bO,\gO,\gO,\gO,\gO,\gO,\bO,\gO,\gO,\gO,\gO,\gO,\gO,\rO,\gO,\gO,\gO,\gO,\gO),\\
\u_9 \:\, = (\gO,\gO,\gO,\gO,\gO,\gO,\bO,\gO,\gO,\gO,\gO,\gO,\bO,\gO,\gO,\gO,\gO,\gO,\gO,\bO),  	&\v_9\:\,= (\gO,\bO,\gO,\gO,\gO,\gO,\gO,\gO,\rO,\gO,\gO,\gO,\gO,\gO,\gO,\bO,\gO,\gO,\gO,\gO),\\
\u_{10}=(\gO,\gO,\gO,\gO,\gO,\gO,\bO,\gO,\gO,\gO,\gO,\gO,\gO,\bO,\gO,\gO,\gO,\gO,\gO,\bO), &\v_{10}=(\gO,\gO,\rO,\gO,\gO,\gO,\gO,\gO,\gO,\bO,\gO,\gO,\gO,\gO,\gO,\bO,\gO,\gO,\gO,\gO),\\
\u_{11}=(\gO,\gO,\gO,\gO,\gO,\gO,\gO,\rO,\gO,\gO,\gO,\gO,\gO,\bO,\gO,\gO,\gO,\gO,\gO,\gO).  &\v_{11}=(\gO,\gO,\gO,\bO,\gO,\gO,\gO,\gO,\gO,\bO,\gO,\gO,\gO,\gO,\gO,\gO,\bO,\gO,\gO,\gO).\\
\end{array}\]}
\end{table}

Notice that $\u_{11}$ and $\v_7$ are the only two vectors that have two 1's and all the rest include three 1's. The first component is always 0. There is a clear pattern in the perturbations and how the motion of the perturbation is from one vector to the following. We have indicated the observed pattern of 1 or 2 repetitions of the perturbations using red and blue (resp.) in Table~\ref{TablE_perturbation_vectors}. Both in the 1 and 3 mod 4 cases, there are 6 red 1's which correspond to a letter's extra appearance as perturbance only a single time in the period. The blue 1's in the components indicate that the corresponding letter occurs as a part of the perturbation twice. After these observations, the periodicity of this perturbation becomes clearer in vector form. 

We claim the following:
\begin{conjecture} For any $i>j\geq 1$, the $i$-th and $j$-th row that appears in Table~\ref{Table_words} and beyond can be associated with the vectors $\b + \u_{(i-1 \text{ (mod 11)})}$ and $\b + \v_{(j\text{ (mod 11)})}$ for 1 and 3 modulo 4 words, respectively. 
\end{conjecture}
Moreover, we claim, with these two essential periods for the perturbation, that we can locate the maximum absolute coefficients' locations for any $(q;q)_n$ with odd $n$. 

The base word associated with $\b$ has 72 letters, each letter representing 19 consecutive $E_n$ values, where each of these values are an additional $4$ index $n$ shifts of $(q;q)_n$. There are 11 base words corresponding to both odd residue classes mod 4. This translates to shifts of the $n$ of $(q;q)_n$ by $72\times19\times4\times11 = $60{,}192. On top of that, the total amount of perturbation is 32 letters, reflecting another 2{,}432 to the shifts of the subindex of $(q;q)_n$, respectively. All together, the period of $E_n$'s is conjectured to be 62{,}624. Lastly, the sum of $E_n$'s is 19{,}787. Now we can start stating our conjectures for the maximum absolute coefficients' location, using~\eqref{Location}.

\begin{conjecture}\label{Conj_Small}We have,\begin{align*}
L(5909 + 68324k) &= \frac{(5909 + 62624k)(5910 + 62624k)}{4} - \frac{3735}{2} - 19787k ,\\
L(391 + 62624k) &= \frac{(391 + 62624k)(392 + 62624k)}{4} - 124 - 19787 k ,
\end{align*} for any $k\geq 0 $.
\end{conjecture}

Note that $D_{5909}=3735/2$ and $D_{391}=124$. Although these formulas might seem limited, one can easily use these explicit values with $E_n$'s as represented in Table~\ref{Table_words} to find the location of the maximum coefficient of any odd value greater than 5909. The location calculations for the $n\leq 5000$ already exist in the Online Encyclopedia of Integer Sequences, under the sequence number~A160089~\cite{OEIS}, and one can easily recover the gap between 5000 and 5909, by our above mentioned calculations. 

We list similar conjectures that correspond to the start of each word of the periods of 1 and 3 modulo 4 cases:

\begin{conjecture}\label{Conj_Big} For any $k\geq 0$, $r=0,1,2,\dots 10$, let \begin{equation*}A(k,r) = 5700r + 62624k - 76\, \delta_{r,11},\text{  and\hspace{.2cm}}B(k,r) = 5700r + 62624k - 76\, \delta_{r\geq 7},\end{equation*} where $\delta_{i,j}$ is the Kronecker delta and $\delta_{i\geq j}$ yields 1 if $i\geq j$, and 0, otherwise. Then,
\begin{align*}
L(5909 + A(k,r)) &= \frac{(5909 + A(k,r) )(5910 + A(k,r) )}{4} - D_{5909 + A(0,r)} - 19787k,\\
L(391 + B(k,r)) &= \frac{(391 + B(k,r))(392 + B(k,r) )}{4} - D_{391 + B(0,s)} - 19787k,\\
\end{align*}where the full list of needed seed values of $D_n$'s are
\[\begin{array}{cc|cc}
D_{5909}=1867.5 		& 	D_{11609} 	= 3668.5,	&  	D_{391} =124,& 	D_{6091}=1925,	\\
D_{17309}=	5469.5,		&	D_{23009}	=7270.5	,	&	D_{11791}=3726 ,&	D_{17491}=	5527,\\
D_{28709}	=9071.5	,	&	D_{34409}	=9675.5	,	&	D_{23191}=	7328,&	D_{28891}= 9129,	\\
D_{40109}	=12673.5,	&	D_{45809}	=14474.5,	&	D_{34591}=	10930,&	D_{40215}=12707,	\\
D_{51509}	=16275.5,	&	D_{57209}	= 18076.5,	&	D_{45915}= 14508,	&	D_{51615}=	16309,\\
D_{62833}	=19853.5.	&							&	D_{57315}=	18110.&
\end{array}
\]
\end{conjecture}

In Conjecture~\ref{Conj_Big}, $A(k,r)$ and $B(k,r)$ are used to move to the start of any word. The variable $r$ is used to move within the periods, and $k$ is used for going over full periods. Conjecture~\ref{Conj_Big} reduces to Conjecture~\ref{Conj_Small} for $r=0$. One important note we need to make is the use of the correctional $\delta$ functions. The number 76, as explained after the introduction of the chosen alphabet, corresponds to the number of consecutive $n$'s needed to expand $(q;q)_n$ for forming a single letter. We have already observed that two words in Table~\ref{Table_words} (one for each case) in the period have one less letter than the others. The $\delta$ functions are introduced only to address this issue.

\section{Future directions}

We have handled the odd cases in two groups corresponding to residue classes modulo 4. We can also start with a uniform approach. For all odd $n\geq 35$, let \[\tilde{E}_n = 2 (D_n - D_{n-2}).\] It is clear that one can define a formula for $L(N)$, which would be an analog of \eqref{Location}, using these~$\tilde{E}_n$ values and the~$D_n$'s. The period in this unified approach seems to be larger much larger than 11. Nevertheless, one intriguing claim analogous to Conjecture~\ref{e_size} is:

\begin{conjecture} Let $n$ be an odd number $\geq 61$, then $\tilde{E}_n = 1$ or $3$.
\end{conjecture}

Using the last maximum absolute coefficients we can note that \[\M_{75000}/\M_{74999} = 1.2197054247521000707999488758629989086855723401824\dots.\] Moreover, \[\left(\M_{75000}\right)^{1/75000} = 1.2195493261856946041555028013329081101077213479207\dots,\] which is a closer estimate than Finch's \cite{Finch} approximation for this value, but it neither supports nor disproves Kotesovec's \cite{OEIS} conjecture: \[\lim_{n\rightarrow\infty} \sqrt[n]{\M_n} \approx \lim_{n\rightarrow\infty} \left(\int_0^1 \prod_{j=1}^n 4\sin^2 ( \pi j z) dz \right)^{\frac{1}{2n}} \approx 1.21971547612163368901359933\dots. \]

We also checked the possible fitting functions and their limits using our data to understand these constants better. Starting from $n=55000$, and using increments of 1000, we have fitted the function $\M_{n}/\M_{n-1}$ and $\sqrt[n]{\M_n}$ to understand their aymptotic behaviour better. We assume that these functions asymptotically look like a series of form \[a_0 + \frac{a_1}{n} + \frac{a_2}{n^2}+\dots,\] where the major contribution comes from the first three coefficients. With that the estimates for the limiting constants are as follows:
\begin{align*}
\lim_{n\rightarrow\infty} \M_{n}/\M_{n-1}&\approx 1.2197154446807031938302328106364792485062669165\dots,\\
\lim_{n\rightarrow\infty} \sqrt[n]{\M_n} &\approx 1.2197063024525922480580886636142211906258238263\dots.
\end{align*}
We can compare these limits with the particular results we have calculated for $n=75000$. This experimentally suggests that, we can trust at least 4 and 3 decimal digits of these last limits for the successive ratios and the successive roots, respectively. 

\section{Acknowledgment}

The authors would like to sincerely thank the interest and encouragement of Jakob Ablinger, Stefano Capparelli, Ralf Hemmecke, Veronika Pillwein, Carsten Schneider, Zafeirakis Zafeirakopoulos, and Wadim Zudilin. We are particularly grateful to Christoph Koutschan and Elaine Wong for their interest, careful reading, and helpful comments on the manuscript.

The authors would like to acknowledge Johannes Bl\"{u}mlein and Carsten Schneider for allowing us to use their servers for the initial calculations at the Research Institute for Symbolic Computation (RISC), Ralf Wahner and  Werner Danielczyk-Landerl for their help in the storage of the early calculations at RISC, and Johann Messner from the MACH2 team for his help in carrying out the calculations. 

Research of the first author is partly supported by the Simons foundation, Award ID: 308929. Research of the second author is supported by the Austrian Science Fund FWF, SFB50-07, SFB50-09 and SFB50-11 projects.

\end{document}